\documentclass[12pt]{amsart}
\usepackage{amssymb} 
\usepackage[mathscr]{eucal}
\usepackage{epsf,epsfig}

\numberwithin{equation}{section} 

\makeatletter
\def\LaTeX{\leavevmode L\raise.42ex
    \hbox{\kern-.3em\size{\sf@size}{0pt}\selectfont A}\kern-.15em\TeX}
\makeatother

\newcommand{\BibTeX}{{\rm B\kern-.05em{\sc i\kern-.025emb}\kern-.08em\TeX}}

\newtheorem{thm}{Theorem}
\newtheorem{prop}[thm]{Proposition}

\theoremstyle{definition}
\newtheorem{defn}[thm]{Definition}
\newtheorem{exmp}[thm]{Example}

\newtheorem{rem}[thm]{Remark} 

\newenvironment{ackn}{\medskip \noindent \small
{\sl Acknowledgments.}}{\bigskip}

\makeatletter
\def\alphenumi{%
  \def\theenumi{\alph{enumi}}%
  \def\p@enumi{\theenumi}%
  \def\labelenumi{(\@alph\c@enumi)}}
\makeatother

\newcommand{\marginnote}[1]
{
}
\newcounter{cv}

\newcounter{gm}

\newcounter{bk}

\begin{document}
\bigskip
\hfill{\it To Vladimir Igorevich Arnold}

\hfill{\it on the occasion of his 70th birthday}
\bigskip
\bigskip

\title[Burgers and Euler]{Shock waves for the Burgers
equation and \\curvatures of diffeomorphism groups}
\author{Boris Khesin}
\address{B.K.: Department of Mathematics, University of Toronto, 
ON M5S 2E4, Canada} 
\email{khesin@math.toronto.edu} 
\author{Gerard Misio\l ek}
\address{G.M.: Department of Mathematics, University of Notre Dame, 
IN 46556, USA}
\email{gmisiole@nd.edu} 

\date{February 1, 2007}


\newcommand{\abs}[1]{\lvert#1\rvert}

\newcommand{\blankbox}[2]{%
  \parbox{\columnwidth}{\centering
    \setlength{\fboxsep}{0pt}%
    \fbox{\raisebox{0pt}[#2]{\hspace{#1}}}%
  }%
}

\begin{abstract}
We  establish 
a simple  relation between curvatures of the group 
of volume-preserving diffeomorphisms and the lifespan of potential 
solutions to the inviscid Burgers equation before the appearance
of shocks. 
We show that shock formation corresponds to a focal point
of the group of volume-preserving diffeomorphisms regarded
as a submanifold 
of the full diffeomorphism group and, consequently, to 
a conjugate point along a geodesic 
in the Wasserstein space of densities. This establishes 
an intrinsic connection between ideal Euler
hydrodynamics (via Arnold's approach), shock formation in 
the multidimensional Burgers equation 
and the Wasserstein geometry of the space of densities. 
\end{abstract}

\maketitle

\section*{Introduction} 
In the famous 1966 paper \cite{Arn} V. Arnold described the dynamics of 
an ideal fluid as a geodesic flow on the group of volume-preserving 
diffeomorphisms of a fixed domain equipped with the metric defined by 
the kinetic energy. He also showed how 
sectional curvature of this group enters 
the problem of Lagrangian stability of Eulerian fluid motions. 
In this paper we are concerned with the exterior geometry of the 
group of volume-preserving diffeomorphisms considered as 
an infinite-dimensional submanifold of the group of all diffeomorphisms. 
We study its second fundamental form and introduce the associated 
shape operator. The focal points of this infinite-dimensional submanifold 
turn out to correspond to shocks forming in potential solutions of the 
inviscid Burgers equation on the underlying finite-dimensional domain. 
This provides a common geometric framework for both Eulerian 
hydrodynamics of ideal fluids and the phenomenon of shocks of 
the inviscid Burgers equation.

\medskip

More precisely, let $M$ be a compact $n$-dimensional Riemannian manifold. 
Consider the group $\mathcal{D}^s(M)$
of Sobolev class diffeomorphisms of $M$ along with its subgroup 
$\mathcal{D}_\mu^s(M)$ of diffeomorphisms preserving the Riemannian 
volume form $\mu$. 
(As usual, if $s> n/2 +1$ both groups can be considered as smooth Hilbert 
manifolds.) 
For a curve $\eta(t)$ in $\mathcal{D}^s(M)$ defined on an interval $[0,a]$ 
its $L^2$-energy is given by 
\begin{equation} \label{energy} 
E(\eta) 
= 
\frac{1}{2}\int_{0}^a \|\dot{\eta}(t)\|_{L^2}^2 \, dt\,,
\end{equation}
where the norm $\|X\|_{L^2}^2=\langle X, X \rangle_{L^2}$ 
is defined by the $L^2$-inner product 
\begin{equation}\label{L2-metric}
\langle X, Y \rangle_{L^2} 
= 
\int_M \langle X(x), Y(x) \rangle \, dx 
\end{equation}
on each tangent space $T_\eta\mathcal{D}^s(M)$. 
The corresponding (weak) Riemannian metric on  $\mathcal{D}^s(M)$
is right-invariant when restricted to the subgroup
$\mathcal{D}_\mu^s(M)$ of volume-preserving diffeomorphisms, although it is
not right-invariant on the whole of $\mathcal{D}^s(M)$. 

Arnold \cite{Arn} proved that the Euler equation of 
an ideal incompressible fluid occupying
the manifold $M$, 
$$
\partial_t u+\nabla_u u=-\nabla p 
$$
for a divergence-free field $u$, corresponds to 
the equation of the geodesic flow of the above 
right-invariant metric on the group 
$\mathcal{D}_\mu^s(M)$.
He also showed that Lagrangian instability of such flows, 
regarded as geodesic deviation, 
can be estimated in terms of sectional curvatures of  the group of 
volume-preserving diffeomorphisms $\mathcal{D}_\mu^s(M)$ and 
provided first such curvature estimates of this group.

Below, by regarding the group of volume-preserving 
diffeomorphisms as a subgroup of the group of all diffeomorphisms 
(cf. e.g. \cite{EM}), we describe its sectional curvatures  
by means of the second fundamental form
of the embedding $\mathcal{D}_\mu^s(M)\subset \mathcal{D}^s(M)$
and relate its principal curvatures to the distance from its nearest focal
point. Recall that the distance from a submanifold $N$
to the first focal point along a geodesic in the direction
normal to this submanifold
gives a lower bound on the principal curvature 
radius of $N$ (in the subspace 
containing the normal direction): 
$$
\mathrm{dist}(N, nearest~f\!ocal~point) 
=
\min \big|~curvature~ radius~o\!f~N~ \big|
$$
$$
=1/\max \big|~principal~curvature~o\!f~N~ \big|\,.
$$

Hence a lower bound on the distance to a focal point provides an
upper bound for (principal) curvatures of $N$. 
This motivates our   
study of the Riemannian geometry of the  embedding 
$\mathcal{D}_\mu^s(M)\subset \mathcal{D}^s(M)$
as the following sequence of implications, which we make precise below:

$i)\,$ Geodesics in the full diffeomorphism group $\mathcal{D}^s(M)$ 
with respect to the 
above $L^2$-metric are described by solutions of the inviscid Burgers equation.
The geodesics normal to the submanifold 
$\mathcal{D}_\mu^s(M)\subset \mathcal{D}^s(M)$ of 
volume-preserving diffeomorphisms are given by potential Burgers solutions,
see Figure 1.

$ii)\,$ The first focal point along a (potential) Burgers solution 
determines the moment when the shock wave develops. The geometry of 
the initial profile allows one to precisely estimate the lifespan 
of this solution. 

$iii)\,$ Focal points along normal geodesics 
are in one-to-one correspondence with conjugate points along 
the projection of these geodesics to the space of densities 
$\mathcal{P}(M)$ equipped with the Wasserstein $L^2$-metric.

$iv)\,$ On the other hand, sectional curvatures of 
$\mathcal{D}_\mu^s(M)$ can be 
explicitly computed from the Gauss-Codazzi equations using the second 
fundamental form (or, the corresponding shape operator) of the embedding 
$\mathcal{D}_\mu^s(M)\subset \mathcal{D}^s(M)$.

$v)\,$ The location of the first focal point along a normal geodesic 
can be estimated in terms of the spectral radius of the shape operator 
of $\mathcal{D}^s_\mu(M)$.

\medskip

\begin{figure}
\input epsf
\centerline{\epsffile{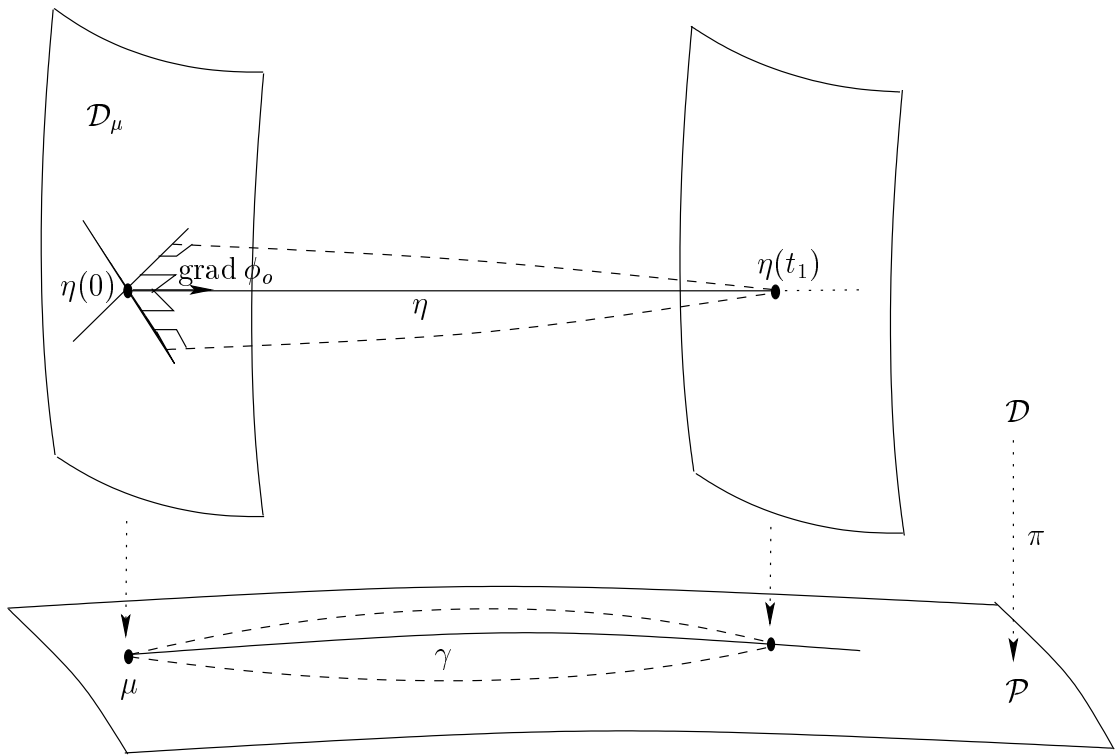}}
\caption{Diffeomorphism group $\mathcal{D}^s(M)$
projects to the space of densities $\mathcal{P}(M)$ with the fiber 
$\mathcal{D}^s_\mu(M)$; normals to $\mathcal{D}^s_\mu(M)$ are horizontal
geodesics, and focal points along them correspond to conjugate points
along geodesics in $\mathcal{P}(M)$.}
\end{figure}


\medskip

In the next three sections we describe the geometry behind  
relations $i)-iii)$, respectively. 
The second fundamental form, the shape operator 
and its spectral properties in
items $iv)-v)$ are described 
in Section \ref{shape}.
A detailed analysis of the shape operator will be 
postponed to a future publication.
In the last section we recall
the properties of asymptotic directions on the diffeomorphism groups
discussed in \cite{KM} within 
the framework of the differential geometry of $\mathcal{D}^s(M)$.
Finally, we mention that 
the connection between the Lagrangian instability of ideal fluids
and negativity of curvatures of the corresponding diffeomorphism group 
was precised by V.Arnold in \cite{Arn, Arn1}, and further explored in 
\cite{l}, \cite{p}.

\bigskip

\section{Burgers potential solutions}

The geometric characterization of the onset of shock waves in terms 
of focal points of the embedding 
$\mathcal{D}_\mu^s(M)\subset \mathcal{D}^s(M)$
is based on the following observation.

\begin{prop}
Geodesics in the group $\mathcal{D}^s(M)$ with respect to 
the $L^2$-metric (\ref{L2-metric}) correspond to solutions of 
the Burgers equation:
each particle moves with constant velocity along a geodesic in $M$.
Geodesics normal to the submanifold $\mathcal{D}^s_\mu(M)$ 
have potential initial conditions.  
\end{prop}

\begin{proof}
First assume for simplicity that $M$ is equipped with a flat metric.
Then the statement that particles move with constant velocity along their 
own geodesics in $M$ follows from the ``flatness'' of the 
$L^2$-metric (\ref{L2-metric}).
To see that the corresponding velocity field satisfies 
the Burgers equation denote the flow of fluid particles by 
$(t,x)\to \eta(t,x)$ and let $u$ be its velocity field 
$$
\frac{d\eta}{dt}(t,x) = u(t,\eta(t,x)), \qquad \eta(0,x)=x. 
$$ 
The chain rule immediately gives 
$$
\frac{d^2\eta}{dt^2}(t,x) 
= 
\big( \partial_t u+ Du\cdot u \big)( t,\eta(t,x))
$$  
and hence the Burgers equation on $M$ 
$$
\partial_t u+ Du\cdot u=0
$$
is equivalent to 
$$
\frac{d^2\eta}{dt^2}(t,x) =0\,, 
$$ 
the equation of freely flying non-interacting particles in $M$.

In the general case, any Riemannian metric on $M$ induces
a unique Levi-Civita $L^2$-connection $\bar{\nabla}$
on $\mathcal{D}^s(M)$, which is determined pointwise 
by the Riemannian connection $\nabla$ on the manifold $M$ itself
(see for example \cite{EM} or \cite{m}). 
Then the same chain rule leads to the Burgers equation
$$
\partial_t u+\nabla_u u=0 
$$
equivalent to the Riemannian version of freely flying
particles: $\bar{\nabla}_{\dot{\eta}}\dot{\eta}=0$.
The latter equation has the explicit form 
$\ddot{\eta} 
+ 
\sum_{ijk}
\Gamma^i_{jk}(\eta) \dot{\eta}^j\dot{\eta}^k
\frac{\partial}{\partial x^i}
=0$, where $\Gamma^i_{jk}$ are the Christoffel symbols of $\nabla$ 
in a local coordinate system $x^1, \dots, x^n$ on $M$. 

Finally, observe that the tangent space to the subgroup 
$\mathcal{D}^s_\mu(M)$ at the identity consists of divergence-free 
vector fields and hence the space of normals is given by gradients 
of $H^{s+1}$ functions on $M$.
Thus horizontal geodesics are the ones whose initial velocities 
are gradient fields $u|_{t=0}={\mathrm{grad}\,\phi_o}$. 
\end{proof}
\bigskip


\section{Burgers shock waves}

The following result describes a connection between focal points and 
formation of shocks for the Burgers equation.

\begin{thm}
The first focal point to the submanifold 
$\mathcal{D}_\mu^s(M)\subset\mathcal{D}^s(M)$ in the direction 
${\mathrm{grad}\,\phi_o}$ is given by the  moment of
the shock wave formation for the Burgers equation with this initial condition 
$u|_{t=0}={\mathrm{grad}\,\phi_o}$.
\end{thm}

\begin{exmp}
In the 1-dimensional case the shock wave solutions of
the Burgers equation appear from the inflection points of the initial 
velocity profile, see Figure 2.
Such points correspond to $u''|_{t=0}=0$, i.e. to $\phi'''_o=0$ 
for $u=\phi'$.
\end{exmp}

\medskip

\begin{figure}
\input epsf
\centerline{\epsffile{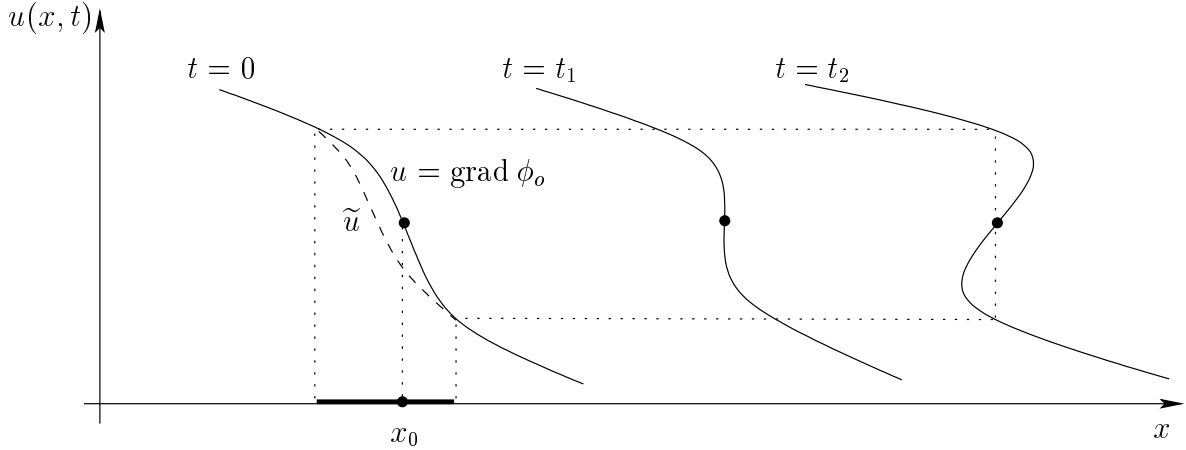}}
\caption{The velocity profile for a solution of the one-dimensional Burgers
equation. An inflection point of the initial profile $u$ generates
the shock wave. The perturbed profile $\tilde u$ leads to a later formation
of the shock wave.}
\end{figure}


\medskip

\begin{proof}
Before the appearence of a shock wave any solution with smooth initial 
data remains smooth and the shock wave is the first moment of non-smoothness.
Furthermore, for a given potential $\phi_o$, any such 
solution is given by a family of diffeomorphisms
$\eta(t): x\mapsto \exp(t\,{\mathrm{grad}\,\phi_o}(x))$ 
parameterized by time $t$. The loss of smoothness occurs 
when the potential $-t\,\phi_o $ ceases to be
$c$-concave with respect to the distance function 
$c(x,y)=d^2(x,y)/2$ 
on the manifold. In this case $\eta(t)$ is no longer a diffeomorphism. 
On the other hand,  as long as the potential is
$c$-concave, the curve $\eta(t)$
remains the shortest curve joining the  diffeomorphisms
$\eta(0)$ and $\eta(t)$, see \cite{McCann}, and focal points on  $\eta(t)$
cannot appear.

Finally, we note that the moment when a shock wave appears corresponds to 
a focal point (i.e., to a caustic), rather than to a cut (or Maxwell) point. 
This follows from the fact that in any
shock wave the first collisions start between neighboring particles,
see e.g. \cite{Arn2, ABB}. To see why such collisions lead to the
existence of a neighboring shorter geodesic, one can employ the
following heuristic argument.
Consider a shock wave generated by a solution 
to the 1-dimensional Burgers equation with the initial velocity 
$u=\mathrm{grad}\,\phi_o$, see Figure 2. 
Assume that this shock wave emerged at time $t_1$
and consider it at time $t_2>t_1$. Consider a perturbation
$\tilde u$
of the original solution which at the initial moment differs from the initial 
velocity $u$ on a small segment. This segment consists of exactly those 
particles  which are glued together in the shock wave between the times 
$[t_1,t_2]$. Perturbations $\tilde u$, whose velocity profile on 
the $(x,u)$-plane (see Figure 2) lies below that of $u$ 
correspond to smooth perturbations of the initial velocities 
of these particles, slowing them down to ensure that the shock 
wave develops after time $t_1$ but no later than $t_2$.
Note that such a perturbation defines another solution which
coincides  with the original one after time $t_2$ but represents 
a shorter geodesic.
Indeed, since we decreased the velocity vectors of the particles, 
the length of the new geodesic, being the integral of the velocity squares, 
can only get smaller. By arbitrarily shortening the segment $[t_1,t_2]$
this implies that the first moment $t_1$ of the original shock
wave for the initial condition $u$
indeed defines a focal point, rather than a cut point, since there exists 
a shorter geodesic nearby.
\end{proof}

\begin{rem}
In higher dimensions the shock waves 
are first generated from the special points
of the initial potential $u|_{t=0}={\mathrm{grad}\,\phi_o}\,$: 
singularities of type $A_3$ modulo certain linear 
and quadratic terms in local orthogonal coordinate charts, see \cite{Boga}.
The list of initial singularities, possible bifurcations of the shock waves
and other related questions for the inviscid Burgers equations, can be found in
\cite{Arn2, ABB, FB}. 
A given shock wave  for the Burgers equation  forms at 
a caustic corresponding to a deeper degeneration of this potential, see
e.g. \cite{Arn2}.

The most degenerate caustic points correspond to the endpoints of the spectrum
for the shape operator $Sh_{\mathrm{grad}\phi_o}$, which we
describe below. The  lifespan of a solution before the shock formation
provides an estimate
for the spectrum of this operator: a {\it lower bound} 
for the lifespan, and hence for the distance to the first focal point 
gives the {\it upper bound} for the (principal) curvatures of the submanifold
$\mathcal{D}_\mu^s(M)\subset \mathcal{D}^s(M)$.
\end{rem}

\medskip

\section{Focal points in diffeomorphisms and conjugate points in densities}

Next we would like to show that  conjugate points along geodesics 
in the Wasserstein space of densities are
in 1-1  correspondence with focal points on the group of diffeomorphisms.
Namely, consider the fibration  
$\pi : \mathcal{D}^s(M) \to \mathcal{P}(M)$, where $\pi$ is 
the projection of diffeomorphisms onto the space of (normalized) 
densities on $M$. Thus, 
two diffeomorphisms $\eta_1$ and $\eta_2$ belong to the same fiber 
if and only if   
$\eta_1=\eta_2\circ\varphi$ for some volume-preserving 
diffeomorphism $\varphi$.
(Moreover, this projection can be defined for more general maps by tracing how
they transport the given standard density form on $M$.)
The map $\pi$ is a Riemannian submersion
\cite{Otto} and the above $L^2$-metric on $\mathcal{D}^s(M)$ induces
the Wasserstein (or, Kantorovich--Rubinstein) metric in the space of 
densities $\mathcal{P}(M)$.

The following general observation holds. 
Consider a projection $\pi : \mathcal{D} \to \mathcal{P}$ 
between two (possibly weak) 
Riemannian manifolds which is a Riemannian submersion  
and where $\mathcal{D}$ an infinite-dimensional Lie group. 
Let $\gamma$ be a geodesic  
starting at a point $\gamma_o$ in $\mathcal{P}$ 
and let $\eta$ be 
a horizontal lift of this geodesic to $\mathcal{D}$, 
i.e. a horizontal geodesic in $\mathcal{D}$ 
whose initial point lies in 
the fiber $\pi^{-1}(\gamma_o)$ and whose projection to $\mathcal{P}$ 
is $\gamma$. Assume, in addition, that the metric restricted to 
$\pi^{-1}(\gamma_o)$ is right invariant. 
(In our case,
the fiber $\pi^{-1}(\gamma_o)$ is the group $\mathcal{D}_\mu^s(M)$ of 
volume-preserving diffeomorphisms for $\gamma_o$).

\begin{prop} \label{prop:gamma} 
The points along the geodesic $\gamma\subset \mathcal{P}$ 
conjugate to the 
initial position $\gamma_o$ are in one-to-one correspondence with 
the focal points of the fiber $\mathcal{D}_\mu=\pi^{-1}(\gamma_o)$, 
regarded as a submanifold in $\mathcal{D}$, 
along a horizontal 
geodesic  $\eta$ in $\mathcal{D}$. Moreover, the multiplicities of 
the conjugate points in $\mathcal{P}$ 
coincide with the multiplicities of the corresponding
focal points in $\mathcal{D}$.
\end{prop}

Recall that conjugate points and their generalizations, the focal points, 
arise as singularities of the Riemannian exponential map. 
In a general infinite-dimensional space two types of such points, 
called epi-conjugate (epi-focal) and mono-conjugate (mono-focal) points, 
can be found. 
Roughly speaking, geometric significance of the former has to do 
with covering properties of the exponential map, while the latter 
are responsible for the minimizing properties of geodesics.

\begin{proof}
The statement  follows from the submersion property of $\pi$. 
Since geodesics in the base manifold 
have unique lifts to horizontal geodesics 
in $\mathcal{D}$ 
we can identify the geodesics emanating from $\gamma_o$ that are 
close to $\gamma$ and which come together 
near a conjugate point in $\mathcal P$
with the  horizontal geodesics emanating from $\pi^{-1}(\gamma_o)$,
that are close to $\pi^{-1}(\gamma)$ and come together near 
the focal point in $\mathcal{D}$.

Furthermore, the argument works for both mono- and epi-conjugate 
(or focal) points, since the exponential map is nondegenerate along 
the fibers.
Indeed, shifts along the fibers correspond to right multiplication 
by an element of $\pi^{-1}(\gamma_o)$, 
while the metric in $\mathcal{D}$ is $\pi^{-1}(\gamma_o)$-invariant.
\end{proof}

\begin{rem}
Geodesics in the Wasserstein space $\mathcal{P}(M)$ of 
(smooth) densities are projections 
of the horizontal geodesics in the diffeomorphism group  
$\mathcal{D}^s(M)$. For a given potential $-\phi_o$ 
they are the diffeomorphisms  
$\eta(t): x\mapsto \exp(t\, {\mathrm{grad}\,\phi_o}(x))$ 
parameterized by $t$.
According to the theorem on polar decomposition on manifolds
\cite{McCann} every non-degenerate map $\eta\in \mathcal{D}^s(M)$ 
has a unique decomposition
$\eta=gr\circ vp$ into a ``gradient map''
$gr(x):=\exp_x({\mathrm{grad}\phi_o})$ for a $c$-concave potential $-\phi_o$ 
and a volume-preserving map $vp$. 
Consequently, the projection of the family of diffeomorphisms $\eta(t): 
x\mapsto \exp(t\, {\mathrm{grad}\phi_o}(x))$ remains 
a shortest geodesic in  $\mathcal{P}$ (in the space of non-degenerate maps)
as long as the potential $-t\,\phi_o$ remains $c$-concave.
\end{rem}

\bigskip


\section{The shape operator and focal points} \label{shape}

We next turn to the shape operator of the embedding 
$\mathcal{D}_\mu^s(M)\subset \mathcal{D}^s(M)$. 
The eigenvalues of this operator are traditionally called 
the principal curvatures and the corresponding eigenvectors 
of unit norm are the principal directions. 
In our case, at each point in $\mathcal{D}^s_\mu(M)$ there is 
a family of shape operators parametrized by vectors normal 
to the submanifold. 

\begin{defn} 
The {\it shape operator} of the submanifold $\mathcal{D}_\mu^s(M)$ 
at the identity is the operator 
$Sh_{\mathrm{grad}\psi}$ on $T_\mathrm{id}\mathcal{D}^s_\mu(M)$ 
defined by 
$$
\langle Sh_{\mathrm{grad}\psi}(w),v \rangle_{L^2} 
= 
\langle 
S_{\mathrm{id}}(w,v), 
\mathrm{grad}\psi 
\rangle_{L^2},
$$
where $\psi\in H^{s+1}(M)$, 
$w$ and $v$ are in $T_\mathrm{id}\mathcal{D}^s_\mu(M)$ 
and $S_\mathrm{id}$ is the second fundamental form of 
$\mathcal{D}^s_\mu(M)$. 
By right invariance we similarly define the shape operator at any point $\eta$ 
in $\mathcal{D}^s_\mu(M)$. 
\end{defn}

\begin{rem}
Recall that the (weak) Riemannian metric on  $\mathcal{D}^s(M)$
induces on both the group
$\mathcal{D}^s(M)$ and its subgroup $\mathcal{D}^s_\mu(M)$ the
unique Levi-Civita connections 
$\bar{\nabla}$ and $\tilde{\nabla}=P_\eta\bar{\nabla}$, 
where 
$P_\eta := R_\eta\circ P_{\mathrm{id}}\circ R_{\eta^{-1}}$ 
is the Hodge projection $P_\mathrm{id}$ onto divergence-free 
vector fields on $M$ conjugated with the right translation 
by $\eta$ (see \cite{EM} for more details).

The second fundamental form $S$ of the submanifold $\mathcal{D}^s_\mu(M)
\subset \mathcal{D}^s(M)$ is the difference of the two connections,
which at the identity is given by 
$$
S_\mathrm{id}(w,v) = Q_\mathrm{id}\nabla_w v \,,
$$
where 
$
Q_\mathrm{id}(w) := w - P_\mathrm{id}(w) 
= 
\mathrm{grad}(\Delta^{-1}\mathrm{div}\, w) 
$ 
is the Hodge projection onto the gradient fields on $M$. 
The second fundamental form, as well as the curvature tensors 
of the connections on 
$\mathcal{D}^s(M)$ and $\mathcal{D}^s_\mu(M)$, 
are bounded multilinear operators satisfying 
the Gauss-Codazzi equations, 
see \cite{m}. 
\end{rem}

The above gives an explicit formula for the shape operator at the identity 
$$
Sh_{\mathrm{grad}\psi}(w) 
= 
- P_{\mathrm{id}} \nabla_w \mathrm{grad}\psi \,.
$$
One can check that the shape operator is bounded and self-adjoint.

The next result provides an estimate of the distance to the first focal 
point along a normal geodesic in terms of the spectral radius of 
$Sh_{\mathrm{grad}\phi_o}$.

\begin{thm} \label{ThmF} 
Suppose that $M$ is a compact manifold of non-negative 
sectional curvature. 
Let $\eta(t)$ be the geodesic normal to $\mathcal{D}^s_\mu(M)$ 
at $\eta_o$ with 
$\dot{\eta}(0)=\mathrm{grad}\phi_o \circ\eta_o$ 
and let $\lambda = \|Sh_{\mathrm{grad}\phi_o}\|_{L(\mathcal{H})}$. 
Then $\eta(t)$ contains a focal point 
in the interval $0<t\leq 1/\lambda$. 
\end{thm} 

In this case our strategy to show existence of focal points along 
geodesics normal to the fibre $\mathcal{D}^s_\mu(M)$ 
is to obtain estimates on
the second variation of the $L^2$-energy (\ref{energy}). 
The latter is given by 
the following formula 
$$
E''(\eta)(W,W)
=
\int_0^a \Big(
\|\bar{\nabla}_{\dot{\eta}}W\|_{L^2}^2
-
\langle \bar{R}(W,\dot{\eta})\dot{\eta},W \rangle_{L^2}
\Big) dt
-
\langle Sh_{\mathrm{grad}\phi_o}(w_o), w_o \rangle_{L^2} 
$$
for any vector field $W$ along $\eta$ such that $W(0)=w_o\circ\eta$ is tangent 
to the fibre $\mathcal{D}^s_\mu(M)$ and $W(a)=0$.
Here $ \bar{R}$ is the curvature of
the connection $\bar{\nabla}$ on the group
$\mathcal{D}^s(M)$. We discuss this theorem and its proof in detail in a 
future publication.

We point out that the conclusion of Theorem \ref{ThmF} holds 
provided that the geodesic $\eta$ 
is defined on the interval 
$[0,1/\lambda]$. 
Furthermore, this geodesic
does not minimize the $L^2$-distance between 
$\mathcal{D}^s_\mu(M)$ and any $\eta(t)$ with $t>1/\lambda$ 
and hence the focal point is in fact mono-focal.

Another corollary of the same formula for the second variation is 
as follows.

\begin{prop}\label{cor:non-pos}
Suppose now that $M$ has non-positive sectional curvature. 
Consider any initial vector $\mathrm{grad}\phi_o$ with
a convex function $\phi_o$. Then 
there can be no focal points on the corresponding geodesic.
\end{prop}

Indeed, for such a manifold $M$ and  
any $W$ one obtains 
$$ 
E''(\eta)(W,W)
\geq 
\int_0^a \|\bar{\nabla}_{\dot{\eta}}W\|_{L^2}^2 \,dt
+
\int_{M} \langle \nabla_{w_o}\mathrm{grad}\,\phi_o, w_o \rangle \,.
$$ 
In particular, if $\phi_o$ is convex then its Hessian 
$
\langle \nabla_{w_o}\mathrm{grad}\,\phi_o, w_o \rangle
$ 
is non-negative and so $E''(\eta)(W,W) \geq 0$, implying
the absence of focal points between $\eta(0)$ and any other $\eta(a)$.

\bigskip
\section{Asymptotic directions} 

Asymptotic directions for the group of volume-preserving diffeomorhisms
also appear naturally in the context of the exterior geometry of this group.
Recall, that a vector tangent to a Riemannian submanifold is
{\it asymptotic} if the (vector-valued) 
second fundamental form evaluated on it is zero. 
The geodesics issued in the direction of this vector, one 
in the submanifold and the other in the ambient manifold, 
have a second order of tangency. (Note that in general two geodesics 
with a common tangent will have only a simple, 
i.e. first order, tangency.) Asymptotic vectors for the submanifold
of volume-preserving diffeomorhisms among all diffeomorhisms are 
given at the tangent space 
to the identity by vector fields $X$
satisfying $\mathrm{div}\,\nabla_{X}X = \mathrm{div}\,X = 0$, see 
\cite{br, KM}. Such fields rarely exist. For instance, 
for 
diffeomorhism groups of surfaces one has the following 
sufficient condition for non-existence of such directions. 

\begin{thm}\cite{KM}
If $M$ is a compact closed surface of nowhere zero curvature, 
then $\mathcal{D}_\mu(M)$ has no asymptotic directions. 
\end{thm}

\begin{proof} Consider the square length function $f:= g(X,X)$ on $M$. At its 
maximum point $x_0$ we get from $df(x_0)=0$ that 
the Jacobi matrix $DX$ is degenerate at $x_0$. This implies that 
$\mathrm{tr}(D X)^2(x_0)=-2 \det{[DX(x_0)]}=0$ for a divergence-free field $X$.

Now employing the identity
$$\mathrm{div}\nabla_X X
=
r(X,X) + \mathrm{tr}(D X)^2\,,
$$
where $r$ stands for the Ricci curvature of the metric $g$, 
which holds for any divergence-free 
vector field $X$ on $M$, we find that 
$$
\mathrm{div}\, \nabla_X X(x_0)
= K(x_0) \, g(X, X)(x_0)\,,
$$
since in two dimensions the Ricci curvature $r$ and the Gaussian 
curvature $K$ coincide.
However, for an asymptotic field $X$ this implies that
$$
0 = K(x_0) \, g(X,X)(x_0)
$$
contradicting the assumption $K\neq 0$ on $M$, and in particular
at the point $x_0$.
\end{proof}

\begin{rem}
Vanishing of the second 
fundamental form on asymptotic vectors implies vanishing of the projections 
of this form to any gradient direction.
In particular, such vectors must be asymptotic simultaneously for 
all the projections of the second fundamental form described above.
This implies the following sufficient condition for the
non-existence of asymptotic directions: if at least in one of the gradient 
projections the second 
fundamental form is sign-definite, then there are no asymptotic directions. 
The latter sufficient condition can also 
be given in terms of gradient solutions
of the Burgers equation: in such a gradient direction all shocks develop
only as $t$ changes in the one direction (say, increases), 
and do not develop as $t$ changes in the other direction. 
\end{rem}


\begin{ackn}
We thank C.\,Villani for fruitful discussions and  
R.\,Wendt for drawing figures.
B.K. is grateful to the MSRI at Berkeley for kind hospitality.  
This research was partially supported by an NSERC  grant
and conducted during the period the first 
author was employed by the Clay Mathematics Institute as a 
Clay Book Fellow. 
\end{ackn}


\bibliographystyle{amsplain}

\end{document}